\documentclass[english]{article}
\usepackage[T1]{fontenc}
\usepackage[latin9]{inputenc}
\usepackage{geometry}
\usepackage{xcolor}
\usepackage{amsmath}
\usepackage{amsthm}
\usepackage{amssymb}
\usepackage{setspace}
\makeatletter
\theoremstyle{plain}
\newtheorem{thm}{\protect\theoremname}
\theoremstyle{plain}
\newtheorem{question}[thm]{\protect\questionname}
\theoremstyle{plain}
\newtheorem{prop}[thm]{\protect\propositionname}
\theoremstyle{remark}
\newtheorem{rem}[thm]{\protect\remarkname}
\theoremstyle{definition}
\newtheorem{defn}[thm]{\protect\definitionname}
\theoremstyle{plain}
\newtheorem{lem}[thm]{\protect\lemmaname}
\theoremstyle{remark}
\newtheorem{conclusion}[thm]{\protect\conclusionname}

\makeatother

\usepackage{babel}
\providecommand{\conclusionname}{Conclusion}
\providecommand{\definitionname}{Definition}
\providecommand{\lemmaname}{Lemma}
\providecommand{\propositionname}{Proposition}
\providecommand{\questionname}{Question}
\providecommand{\remarkname}{Remark}
\providecommand{\theoremname}{Theorem}

\begin{document}
\title{A Bubble-breaking Phenomenon in the Variation of a Swarm Communication Network}
\author{Bohuan Lin, Jian Gao}
\maketitle
\begin{abstract}
We discuss a specific circumstance in which the topology of the communication
network of a robotic swarm has to change during the movement. The
variation is caused by a topological obstruction which emerges from
certain geometric restrictions on both the environment and the swarm.
\end{abstract}

\section{Introduction}

\subsection{Networks and system dynamics}

Communication among individuals in a group of agents (fish/birds/robots/people)
is the foundation for forming system behaviours. The ``communication''
here can be any exchange/flow of information in any form. In flocking/schooling
of birds/fish it may refer to the sensing (by the group members) of
changes in positions/velocities of nearby companions \cite{hemelrijk2012schools}.
In a multirobot system it can mean the transmission of signals among
the robots. In epidemic spreading it is then the transfer of viruses/bacteria. 

Representing with vertices the agents and with edges the established
communications yields a communication network of the group. These
networks serve as basic mathematical structures upon which system
dynamics are built \cite{strogatz2001exploring,olfati2006flocking,holme2012temporal,pastor2015epidemic}.
In nature as well as in practise, the ability of an agent to set up
communication with other individuals is usually limited. For the case
of interest in this note, it is the maximal distance to set up message
channels between two robots. As a consequence, the motion of the group
may in turn cause changes to the topology of the network.

\subsection{A bubble-breaking phenomenon}

Put an iron ring (rigid) into a balloon (not inflated) and seal the
valve. No matter how hard you try, it is impossible to stretch and
attach the entire balloon onto the ring without breaking it. In other
words, if you force as large a part of the balloon as possible to
be attached onto the ring, then the rest part of the balloon will
be terribly stretched and the continuum of the material may eventually
break down. We call this issue a $bubble-breaking$ phenomenon, and
its mathematical essense can have implications on various problems
from different contexts. For example, in stability theory it implies
the fact that, on any sphere enclosing an asymptotically stable limit
cycle $\mathcal{P}$ of a (smooth) flow $\varphi$ in $\mathbb{R}^{3}$,
there always exists some point $p$ such that the trajectory $t\mapsto\varphi^{t}(p)$
does not converge to $\mathcal{P}$. We encoutered this fact in \cite{yao2022topological}
when studying path-following control. In this note, we demonstrate
a case in which the topology of the communication network of a robotic
swarm $has$ $to$ change during its movement. As we will see, the
mathematics behind also reflects the essense of the bubble-breaking
phenomenon. 

The description of the case is detailed in Subsection \ref{subsec:Description-of-the-question}
where the central problem of this study is stated as Question \ref{que:main question},
and then Proposition \ref{prop:main result} is proved in Subsection
\ref{subsec:Answer-to-the question} as answer to the question. The
analysis in Section \ref{sec:The case for Communication-Networks}
assumes a special structure from the initial positions of the robots.
In Section \ref{sec:Further-Discussion} we show a natural condition
on the intitial positions under which the analysis can be applied,
and the main result in this part is demonstrated as Proposition \ref{prop:natural condition and triangulation}.
In the end, by combining Propositions \ref{prop:main result} and
\ref{prop:natural condition and triangulation} we draw the final
conclusion as Theorem \ref{thm:network-variation under natural-condition}.

\section{\label{sec:The case for Communication-Networks}A Case for the Communication
Network of a Robotic Swarm}

\subsection{\label{subsec:Description-of-the-question}Description of the question}

Imagine that there is a facility $\mathrm{\mathfrak{P}}$ (floating
in the space) which occupies an area of a solid torus
\[
\mathrm{\mathfrak{P}}=\big\{(r\cos\theta,r\sin\theta,z)\big|\theta\in[0,2\pi],\,z^{2}+(r-\frac{1}{10})^{2}\leq\epsilon^{2}\big\}
\]
with $\epsilon<\frac{1}{100}$. Note that the circle
\[
\mathrm{\mathcal{P}}=\{x^{2}+y^{2}=\frac{1}{10},z=0\}
\]
is the central axis of the solid torus $\mathrm{\mathfrak{P}}$. Suppose
that the robots can only move in the area out of $\mathrm{\mathfrak{P}}$,
and each robot can only communicate with those within the distance
$\delta<\epsilon$. 

A movement of a group of $n$ robots can be represented as a continuous
map
\[
\mathcal{R}:\mathbb{R}\ni t\mapsto\big(r_{1}(t),...,r_{n}(t)\big)\in(\mathbb{R}^{3})^{n}
\]
with $r_{i}(t)$ being the position of the $i$th robot at the moment
$t$. The communication network among the individuals of the group
can be represented by an $n\times n$ ``$0$-$1$'' matrix $[p_{ij}]$.
To be precise, if there is a message channel built up between the
$i$ th and the $j$ th robots with $i<j$, then we set $p_{ij}=1$
and otherwise $p_{ij}=0$. Since (the structure of) the network may
change over time, we use the symbol $[p_{ij}]_{t}:=[p_{ij}^{t}]$
to denote the network at the moment $t$.

Denote by $o_{i}$ the initial poision of the $i$th robot, i.e. $o_{i}=r_{i}(0)$.
For simplicity, we \textbf{assume} that the initial communication
network $[p_{ij}]_{0}$ induces a triangulation of $S^{2}$. To be
precise, let $a_{ij}$ be the minimal geodesic connecting $o_{i}$
and $o_{j}$ on $S^{2}$. If $p_{ik}^{0}=p_{ij}^{0}=p_{jk}^{0}=1$,
the arcs $a_{ij},a_{jk},a_{ik}$ enclose a geodesic triangle $\sigma_{ijk}$
(including the interior). If two different triangles $\sigma$ and
$\sigma'$ have nonempty intersection, their overlap is either an
edge $a$ or a vertex $o$. Suppose that at the moment $\bar{t}>0$,
the robots land on $\mathrm{\mathfrak{P}}$, i.e. $r_{i}(\bar{t})\in\mathrm{\mathfrak{P}}$.
We consider the question that whether the topology of the network
changes during $this$ movement.
\begin{question}
\label{que:main question}$[p_{ij}]_{t}\equiv[p_{ij}]_{0}$ for $t\in[0,\bar{t}]$?
\end{question}

The answer is not surprising: it has to change. An explanation is
given in the next subsection.

\subsection{\label{subsec:Answer-to-the question}Answer to the question}

Associated to the matrix $[p_{ij}]_{0}$ there is a $2$ dimensional
simplicial complex in $\mathbb{R}^{3}$. For $i,j,k$ such that $p_{ik}^{0}=p_{ij}^{0}=p_{jk}^{0}=1$,
let $\Delta_{ijk}$ be the triangle in $\mathbb{R}^{3}$ with vertices
$o_{i}$, $o_{j}$ and $o_{k}$. $\Delta_{ijk}$ can be seen as a
linear approximation to the curved triangle $\sigma_{ijk}$ on $S^{2}$,
and 
\[
K=\bigcup\Delta_{ijk}
\]
is a polytope with a triagulization given by $\{\Delta_{ijk}\}$.
Moreover, $\Delta_{ijk}$ is homeomorphic to $\sigma_{ijk}$ via the
radial projection. That is, for each $p$ on $\Delta_{ijk}$, there
exists a unique $b_{p}\geq1$ such that $b_{p}\cdot p\in\sigma_{ijk}$,
and the map 
\[
\Delta_{ijk}\ni p\xrightarrow{\tau_{ijk}}b_{p}\cdot p\in\sigma_{ijk}
\]
is a homeomorphism. Piecing together the maps $\tau_{ijk}$ we get
a homeomorphism $\tau$ from $K$ to $S^{2}$ with $\tau\big|_{\Delta_{ijk}}=\tau_{ijk}$.
With small $\delta$, $K$ also encloses $\mathcal{P}$. Then 
\begin{equation}
h_{t}(w)=(1-t)\cdot w+t\cdot\tau^{-1}(w)\label{eq:h_t}
\end{equation}
is a continuous map from $S^{2}\times[0,1]$ to $\mathbb{R}^{3}-\mathcal{P}$
and it is a homotopy between the inclusion $\iota_{S^{2}}=h_{0}$
and $\tau^{-1}=h_{1}$.

The movement $\mathcal{R}$ induces a homotopy 
\begin{equation}
\mathcal{K}:K\times[0,1]\rightarrow\mathbb{R}^{3}\label{eq:K_t}
\end{equation}
by sending each point $p=s_{i}o_{i}+s_{j}o_{i}+s_{k}o_{i}$ in $\Delta_{ijk}$
to the point $s_{i}r_{i}(t)+s_{j}r_{j}(t)+s_{k}r_{k}(t)$. Here $s_{i},s_{j},s_{k}\geq0$
are the weights with $s_{i}+s_{j}+s_{k}=1$, and $\mathcal{K}$ is
well-defined since $\{\Delta_{ijk}\}$ is a triangulation of $K$.
Note that $\Delta_{ijk}^{t}=\mathcal{K}_{t}(\Delta_{ijk})$ is also
a ``triangle'': it is the convex hull of the points $r_{i}(t)$, $r_{j}(t)$
and $r_{k}(t)$ in $\mathbb{R}^{3}$, and its diameter equals to the
largest distance between these points.

We conclude the discussion with the following proposition.
\begin{prop}
\label{prop:main result}There exists $t\in[0,\bar{t}]$ such that
$[p_{ij}]_{t}\neq[p_{ij}]_{0}$.
\end{prop}

\begin{proof}
We argue by contradiction. Assume that $[p_{ij}]_{t}=[p_{ij}]_{0}$
for all $t\in[0,\bar{t}]$. It means that, if $p_{ij}^{0}=1$, then
the distance between $r_{i}(t)$ and $r_{j}(t)$ is always no larger
than $\delta$. As a consequence, the diameter of the ``triangle''
$\Delta_{ijk}^{t}=\mathcal{K}_{t}(\Delta_{ijk})$ is no more than
$\delta$. Since the robots are moving outside $\mathfrak{P}$ and
$\delta$ is smaller than $\epsilon$ (the radius of $\mathfrak{P}$),
we know that $\Delta_{ijk}^{t}$ has no intersection with $\mathcal{P}$.
Since 
\[
\mathcal{K}_{t}(K)=\bigcup\mathcal{K}_{t}(\Delta_{ijk}),
\]

this means that $\mathcal{K}$ is a continuous map from $K\times[0,1]$
to $\mathbb{R}^{3}-\mathcal{P}$. Note that it follows directly from
Eq.(\ref{eq:K_t}) that $\mathcal{K}_{0}$ is the inclusion of $K$
into $\mathbb{R}^{3}-\mathcal{P}$. Then we get a homotopy $\bar{h}$
between $\iota_{S^{2}}$ and $\mathcal{K}_{\bar{t}}\circ\tau^{-1}$
by letting $\bar{h}_{t}=h_{t}$ for $t\in[0,1]$ and $\bar{h}_{t}=\mathcal{K}_{t-1}\circ h_{1}$
for $t\in[1,\bar{t}+1]$. 

Since the diameter of $\Delta_{ijk}^{\bar{t}}$ is no larger than
$\delta$ and $r_{i}(\bar{t})\in\mathfrak{P}$ for all $i\in\{1,...,n\}$,
each $\Delta_{ijk}^{\bar{t}}$ lies in the $\delta$-neighbourhood
$\mathcal{U}$ of $\mathfrak{P}$ and then the image
\[
\bar{h}_{\bar{t}+1}(S^{2})=\bigcup\mathcal{K}_{\bar{t}}(\Delta_{ijk})=\bigcup\Delta_{ijk}^{\bar{t}}
\]
also lies in $\mathcal{U}$. Considering the radii of $\mathcal{P}$
and $\mathfrak{P}$, $\mathcal{U}$ is a thickened $2$-torus lying
in $\mathbb{R}^{3}-\mathcal{P}$, i.e. $\mathcal{U}\cong\mathfrak{P}\times(-\delta,\delta)$
and $\mathcal{U}\cap\mathcal{P}=\emptyset$. Since $S^{2}$ is simply
connected, $\bar{h}_{\bar{t}+1}$ factors through the universal covering
$\mathbb{R}^{3}$ (contractible) of $\mathcal{U}$. The induced homomorphism
$\bar{h}_{(\bar{t}+1),*}$ from $H_{2}(S^{2})$ to $H_{2}(\mathbb{R}^{3}-\mathcal{P})$
then factors as
\[
\bar{h}_{(\bar{t}+1),*}:H_{2}(S^{2})\rightarrow H_{2}(\mathbb{R}^{3})\rightarrow H_{2}(\mathcal{U})\xrightarrow{\text{inclusion}_{*}}H_{2}(\mathbb{R}^{3}-\mathcal{P})
\]

and is therefore a null morphism (since $H_{2}(\mathbb{R}^{3})\cong\{0\}$),
and hence so is $\iota_{S^{2},*}$ due to the homotopy. However, for
any $p_{0}\in\mathcal{P}$, $S^{2}$ is a deformation retract of $\mathbb{R}^{3}-\{p_{0}\}$
and hence 
\[
H_{2}(S^{2})\xrightarrow{\iota_{S^{2},*}}H_{2}(\mathbb{R}^{3}-\mathcal{P})\xrightarrow{\text{inclusion}_{*}}H_{2}(\mathbb{R}^{3}-\{p_{0}\})
\]
is an isomorphism, meaning $\iota_{S^{2},*}$ should not be trivial,
yielding a contradiction.
\end{proof}

\section{\label{sec:Further-Discussion}Further Discussion}

In Section \ref{sec:The case for Communication-Networks} we assume
that the communication network of the $n$ robots at $t=0$ together
with the initial positions $(o_{1},...,o_{n})$ induces a triangulation
of $S^{2}$. In this section we replace this assumption with a more
natural conidtion. This is formalized as Proposition \ref{prop:natural condition and triangulation},
which is proposed in Subsection \ref{subsec:A-natural-condition}
and (eventually) proved in Subsection \ref{subsec:Existence of Triangulations}.
Subsection \ref{subsec:General-positions on S^2} is devoted for a
technical preparation.
\begin{rem}
A $\emph{position of the swarm}$ ($n$ robots) on $S^{2}$ is an
element $(a_{1},...,a_{n})$ in $(S^{2})^{n}$. We say a network $[p_{ij}]$
with a position $(a_{1},...,a_{n})$ induces a triangulation of $S^{2}$
if and only if connecting all those points $a_{i}$ and $a_{j}$ by
geodesics whenever $p_{ij}=1$ gives a triangulation. 
\end{rem}

\begin{rem}
A sub-network/-graph of the network $[p_{ij}]$ can be represented
as a matrix $[p'_{ij}]$ with the same dimension satisfying the relation
\begin{equation}
p'_{ij}=1\implies p_{ij}=1.\label{eq:sub-network relation}
\end{equation}
\end{rem}

\subsection{\label{subsec:A-natural-condition}A natural condition}

If a subgraph fails to keep its structrure then so does the whole
network. Therefore, Proposition \ref{prop:main result} applies as
long as there is a sub-graph/sub-network inducing a triangulation
on $S^{2}$. In fact, Proposition \ref{prop:main result} still holds
even if such a sub-network induces a triangulation only after an admissible
pertubation on the initial position $(o_{1},...,o_{n})$. Here, an
admissible perturbation refers to a position $(o'_{1},...,o'_{n})\in(S^{2})^{n}$
from which the swarm can move to the actual initial position $(o_{1},...,o_{n})$
while keeping the topology of the network unchanged. To be precise,
\begin{defn}
\label{def:admissible perturbation}$(o'_{1},...,o'_{n})$ is a perturbation
of $(o_{1},...,o_{n})$ admissible to a sub-network $[p'_{ij}]$ if
there is a movement
\begin{equation}
\mathcal{\tilde{R}}:[0,1]\xrightarrow{(\tilde{r}_{1},...,\tilde{r}_{n})}(\mathbb{R}^{3})^{n}-\mathfrak{P}\label{eq:admissible displacement}
\end{equation}
with $\mathcal{\tilde{R}}_{0}=(o'_{1},...,o'_{n})$ and $\mathcal{\tilde{R}}_{1}=(o_{1},...,o{}_{n})$
such that whenever $p'_{ij}=1$, $|\tilde{r}_{i}(t)-\tilde{r}_{j}(t)|<\delta$
holds for all $t\in[0,1]$.
\end{defn}

If a triangulation of $S^{2}$ is induced by a sub-network $[p'_{ij}]$
with the position $(o'_{1},...,o'_{n})$ and the movement $\mathcal{\tilde{R}}$
is admissible to $[p'_{ij}]$, then applying Proposition \ref{prop:main result}
to the ``composed'' movement in which the swarm first takes the movement
$\mathcal{\tilde{R}}$ from $(o'_{1},...,o'_{n})$ to $(o_{1},...,o{}_{n})$
and conitnues with the movement $\mathcal{R}$ will then prove that
the structure of $[p'_{ij}]$ has to change during the whole process.
Since it is unchanged in the first movement $\mathcal{\tilde{R}}$,
we again shows that the structure of the network has to change in
the movement $\mathcal{R}$.
\begin{rem}
\label{rem:actual and virtual initial positions}For convenience,
we will call $o_{i}$ and $o'_{i}$ respectively the actual and the
virtual (initial) positions of the $i$th robot. Similarly, $(o_{1},...,o_{n})$
and $(o'_{1},...,o'_{n})$ are respectively the actual and the virtual
(initial) positions of the swarm.
\end{rem}

Based on the discussion above, we will look for a condition which
allows the initial network $[p_{ij}]_{0}$ to have a sub-graph inducing
a triangulation of $S^{2}$ under an admissible perturbation on $(o_{1},...,o_{n})$.
We formalize it as the following proposition.
\begin{prop}
\label{prop:natural condition and triangulation}Suppose that the
initial positions $\{o_{1},...,o_{n}\}$ of the robots constitute
a $\frac{\delta}{6}$-net on $S^{2}$, and any pair of robots will
set up a message channel if the distance between them is smaller than
$\delta$. Then a sub-graph of the communication network $[p_{ij}]_{0}$
induces a triangulation of $S^{2}$ after a perturbation on $(o_{1},...,o_{n})$
which is admissible to the sub-graph.
\end{prop}

The $condition$ of the set $\{o_{1},...,o_{n}\}$ being a $\frac{\delta}{6}$-net
on $S^{2}$ means that every point on $S^{2}$ is at a distance less
than $\frac{\delta}{6}$ from some point (robot) $o_{i}$. Here we
choose the distance to be the (restriction of the) Euclidean metric
from $\mathbb{R}^{3}$ (on $S^{2}$). That is, the distance between
any two points $o$ and $o'$ on $S^{2}$ is measured by the length
of the vector $o-o'$ in $\mathbb{R}^{3}$. With this metric, the
$\frac{\delta}{6}$-neighbourhood $D(o,\frac{\delta}{6})$ of $o$
on $S^{2}$ is simply the intersection of $S^{2}$ with the $3$-dimension
ball $B(o,\frac{\delta}{6})$ in $\mathbb{R}^{3}$ (with center $o$
and radius $\frac{\delta}{6}$). Here both $D(o,\frac{\delta}{6})$
and $B(o,\frac{\delta}{6})$ are taken as open sets in $S^{2}$ and
$\mathbb{R}^{3}$, respectively. The $condition$ is equivalent to
saying the neighbourhoods $D(o_{i},\frac{\delta}{6})$ constitute
an open cover of $S^{2}$, and we consider it to be natural since
it merely gives a description on the density of the robots on $S^{2}$.
Note that this condition may be coarse in the sense that we could
have taken a (much) larger radius than $\frac{\delta}{6}$, or, say,
a (much) smaller density of the robots. However, giving a finer estimation
on how sparse the robots can be (for inducing a triangulation) is
beyond the scope of this note.

\subsection{\label{subsec:General-positions on S^2}General positions of the
swarm on $S^{2}$}

For any $1\leq i,j,k\leq n$, define a function $f_{ijk}$ on $(S^{2})^{n}$
with
\begin{equation}
f_{ijk}(a_{1},...,a_{n}):=\det[a_{i};a_{j};a_{k}].\label{eq:f_=00007Bijk=00007D}
\end{equation}

By saying a general position (of size $n$) on $S^{2}$ we mean an
element $(a_{1},...,a_{n})$ in $(S^{2})^{n}$ such that $f_{ijk}(a_{1},...,a_{n})\neq0$
for all the triples $(i,j,k)$ with $i<j<k$. Note that when the swarm
is in a general position $(a_{1},...,a_{n})$ on $S^{2}$, the convex
hull of any three robots is a triangle in $\mathbb{R}^{3}$, and its
radial projection on $S^{2}$ is a geodesic triangle. In this subsection
we show that given any (initial) position $(o_{1},...,o_{n})$, with
an arbitrarily small perturbation it yields a general position $(o'_{1},...,o'_{n})$.
More precisely, 
\begin{lem}
\label{lem:general positions are generic}The set $\mathfrak{G}$
of all general positions is open and dense in $(S^{2})^{n}$.
\end{lem}

\begin{proof}
Note that $\mathfrak{G}=(S^{2})^{n}-\mathfrak{C}$ with
\[
\mathfrak{C}=\bigcup_{1\leq i<j<k\le n}f_{ijk}^{-1}(0).
\]
We only need to show that for all triple $(i,j,k)$ with $i<j<k$,
the sets $\mathfrak{G}_{ijk}:=\mathfrak{G}-f_{ijk}^{-1}(0)$ are open
and dense in $(S^{2})^{n}$, and then as their finite intersection
$\mathfrak{G}$ is also dense and open in $(S^{2})^{n}$. 

Since $f_{ijk}^{-1}(0)$ is closed in $(S^{2})^{n}$, $\mathfrak{G}_{ijk}$
is open. To see that $\mathfrak{G}_{ijk}$ is dense, we first look
at the subset $S_{ijk}^{2}$ of $f_{ijk}^{-1}(0)$ defined by containing
all the points $(a_{1},...,a_{n})$ with $a_{i}=a_{j}=a_{k}$. It
is straightforward to check that $S_{ijk}^{2}$ is an embedding of
$(S^{2})^{n-2}$ in $(S^{2})^{n}$, and therefore its complement $\tilde{\mathfrak{G}}_{ijk}:=\mathfrak{G}-S_{ijk}^{2}$
is an open and dense subset of $(S^{2})^{n}$ which contains $\mathfrak{G}_{ijk}$. 

It remains to show that $\mathfrak{G}_{ijk}=\tilde{\mathfrak{G}}_{ijk}-f_{ijk}^{-1}(0)$
is dense in $\tilde{\mathfrak{G}}_{ijk}$. For doing this, we will
verify that $0$ is a regular value of the (restricted) function $\tilde{f}_{ijk}:=f_{ijk}\big|_{\tilde{\mathfrak{G}}_{ijk}}$
on $\tilde{\mathfrak{G}}_{ijk}$. This will imply that 
\[
f_{ijk}^{-1}(0)\bigcap\tilde{\mathfrak{G}}_{ijk}=\tilde{f}_{ijk}^{-1}(0)
\]
is an embedded submanifold with codimension $1$ in $\tilde{\mathfrak{G}}_{ijk}$,
and then as its complement $\mathfrak{G}_{ijk}$ is dense in $\tilde{\mathfrak{G}}_{ijk}$.
Suppose that $(a_{1},...,a_{n})$ is from the set $\tilde{f}_{ijk}^{-1}(0)$.
Without loss of generality, we can assume that $a_{i}\neq a_{j}$.
From \ref{eq:f_=00007Bijk=00007D} it holds $\det[a_{i};a_{j};a_{k}]=0$,
which means that the vectors $a_{i}$, $a_{j}$ and $a_{k}$ locate
on a $2$ dimensional vector subspace $V$ of $\mathbb{R}^{3}$. Let
$u$ be a unit vector perpendicular to $V$. Since $a_{i},a_{j},a_{k}\in S^{2}\cap V$,
we know that $u$ is vertical to these three vectors, and then $(a_{k},u)$
is a tangent vector of $S^{2}$ at $a_{k}$. Check that
\[
\frac{d}{dt}f_{ijk}(a_{1},..a_{k}+tu,..a_{n})=\frac{d}{dt}\det[a_{i};a_{j};tu]=\det[a_{i};a_{j};u]\neq0,
\]
which implies that $(a_{1},...,a_{n})$ is a regular point of $f_{ijk}$.
\end{proof}
As a consequence of Lemma \ref{lem:general positions are generic},
we get a general position $(o'_{1},...,o'_{n})$ from any given $\frac{\delta}{6}$-net
$\{o_{1},...,o_{n}\}$ on $S^{2}$ via an arbitrarily small perturbation.
Furthermore, when $(o'_{1},...,o'_{n})$ is sufficiently closed to
$(o_{1},...,o_{n})$, i.e. the displacement
\begin{equation}
\delta o:=\max_{1\leq i\leq n}\{|o_{i}-o'_{i}|\}\label{eq:perburbation error}
\end{equation}
is small enough, $(o'_{1},...,o'_{n})$ is an admissible perturbation
from $(o_{1},...,o_{n})$, and $\{o'_{1},...,o'_{n}\}$ is also a
$\frac{\delta}{6}$-net on $S^{2}$. We give this an explanation in
the following.

Since there are only finitely many pairs of robots with initial conditions
$|o_{i}-o_{j}|<\delta$, the number
\[
l_{o}:=\max\{|o_{i}-o_{j}|\big||o_{i}-o_{j}|<\delta\}
\]
is strictly smaller than $\delta$. If $\delta o<l_{o}$, then the
``straightline'' movement of the swarm (out of $\mathfrak{P}$) defined
by
\[
\tilde{r}_{i}(t):=(1-t)\cdot o'_{i}+t\cdot o_{i}
\]
gives an admissible perturbation. 

Since the neighbourhoods $D(o_{i},\frac{\delta}{6})$ form a finite
cover of $S^{2}$, for a number $\delta'$ slightly smaller than $\delta$,
the neighbourhoods $D(o_{i},\frac{\delta'}{6})$ also cover $S^{2}$.
To see this, we first take a look at the compact set
\[
A_{1}=S^{2}-\bigcup_{j\neq1}D(o_{j},\frac{\delta}{6}).
\]
Since $A_{1}$ has no intersection with $D(o_{j},\frac{\delta}{6})$
for all $j>1$, it is contained in $D(o_{1},\frac{\delta}{6})$. Due
to the compactness of $A_{1}$, with $\delta'_{1}$ slightly but strictly
smaller than $\delta$, $D(o_{1},\frac{\delta'_{1}}{6})$ also contains
$A_{1}$. Therefore, the sets $D(o_{1},\frac{\delta'_{1}}{6})$ and
$D(o_{j},\frac{\delta}{6})$ for $j>1$ constitute an open cover of
$S^{2}$. Now apply this process inductively to the rest $D(o_{i},\frac{\delta}{6})$.
Suppose that we already get an open cover consisting of the sets $U_{i}=D(o_{i},\frac{\delta'_{i}}{6})$
for $i\leq k$ and $U_{i}=D(o_{i},\frac{\delta}{6})$ for $i\geq k+1$.
The compact set
\[
A_{k+1}=S^{2}-\bigcup_{j\neq k+1}U_{j}
\]
is then contained in $U_{k+1}=D(o_{k+1},\frac{\delta}{6})$. Replace
it with $D(o_{k+1},\frac{\delta'_{k+1}}{6})$ and continue until we
eventually get an open cover $D(o_{i},\frac{\delta'_{i}}{6})$ with
$\delta'_{i}<\delta$ for $i=1,...,n$. Take 
\[
\delta'=\max\{\delta'_{i}\}
\]
and then the neighbourhoods $D(o_{i},\frac{\delta'}{6})$ again cover
$S^{2}$ with $\delta'<\delta$, i.e. $\{o_{1},...,o_{n}\}$ is also
a $\frac{\delta'}{6}$-net. 

Now we take $\delta o<\min\{\frac{\delta-\delta'}{6},l_{o}\}$. Since
$\{o_{1},...,o_{n}\}$ is a $\frac{\delta'}{6}$-net, each $p\in S^{2}$
is contained in some $D(o_{i},\frac{\delta'}{6})$, i.e. $|p-o_{i}|<\frac{\delta'}{6}$.
Therefore, it holds that 
\[
|p-o'_{i}|<\frac{\delta'}{6}+\delta o<\frac{\delta}{6},
\]

and hence we conclude that:
\begin{conclusion}
\label{conc:general-positions for a net}Suppose that the intitial
positions $o_{1},...,o_{n}$ on $S^{2}$ form a $\frac{\delta}{6}$-net
on $S^{2}$. Then by admissible perturbation the swarm is in such
a general position $(o'_{1},...,o'_{n})$ that $\{o'_{1},...,o'_{n}\}$
is also a $\frac{\delta}{6}$-net on $S^{2}$.
\end{conclusion}

\subsection{\label{subsec:Existence of Triangulations}Sub-networks inducing
triangulations}

As a $\frac{\delta}{6}$-net, $\{o_{1},..,o_{n}\}$ is also a a $\frac{\delta'}{6}$-net
on $S^{2}$ for some $\delta'<\delta$. According to Conclusion \ref{conc:general-positions for a net},
we can perturb $(o_{1},...,o_{n})$ to a general position $(o'_{1},...,o'_{n})$
such that $\{o'_{1},...,o'_{n}\}$ is a $\frac{\delta'}{6}$-net,
and, the displacement $\delta o$ defined by (\ref{eq:perburbation error})
is smaller than $\frac{\delta-\delta'}{2}$. By the triangle inequality,
if a subgraph with vertices from $\{o'_{1},...,o'_{n}\}$ has all
the edges shorter than $\delta'$, these edges will then remain shorter
than $\delta$ during the swarm moves back to the actual position
$(o_{1},...,o_{n})$ along the straightlines $\overline{o'_{i}o_{i}}$. 

Instead of directly triangulate $S^{2}$ with $\{o'_{1},...,o'_{n}\}$
as the vertices, we will build a triangulation with the vertices from
another set $\mathcal{V}$, and show that the points in $\mathcal{V}$
are sufficiently close to $\{o'_{1},...,o'_{n}\}$. More specifically,
each point $q$ in $\mathcal{V}$ is obaitned from a point $o'$ in
$\{o'_{1},...,o'_{n}\}$ with $|o'-q|<\frac{\delta'}{6}$, and the
mapping $q\mapsto o'$ is injective.

\subsubsection{The set of vertices $\mathcal{V}$}

By rotating $S^{2}$ if necessary, we assume that $o'_{n}$ is the
north pole $(0,0,1)$. For $h\in[-1,1]$, denote by $L_{h}$ the latitude
\[
L_{h}=S^{2}\cap\{(x,y,z)\big|z=h\}.
\]
Note that $L_{0}$ is the equator. For each $e^{i\theta}\in S^{1}$,
denote by $R_{e^{i\theta}}$ the longitude line intersecting with
$L_{0}$ at the point $(e^{i\theta},0)$. For convenience, we will
refer to $h$ and $e^{i\theta}$ the latitude and the longitude, respectively.
Since $(o'_{1},...,o'_{n})$ is a general position, when $i,j\neq n$,
$o'_{i}$ and $o'_{j}$ have different longitude. Choose a sequence
from $(-1,1)$
\[
h_{-k}<h_{-k+1}<...<h_{0}=0<h_{1}...<h_{k}
\]
such that the distance between any two adjacent latitude lines $L_{j}=L_{h_{j}}$
and $L_{j+1}=L_{h_{j+1}}$ is $\frac{\delta'}{3}$. The integer $k$
is taken in such a way that the distance between $L_{k}$ and $o'_{n}=(0,0,1)$
is smaller than $\frac{2\delta'}{3}$ but no less than $\frac{\delta'}{3}$.
Note that by symmetry this means that the distance between $L_{-k}$
and $(0,0,-1)$ is also between $\frac{\delta'}{3}$ and $\frac{2\delta'}{3}$. 

For each $j\in\{-k,...,k\}$, let $\mathcal{C}_{j}$ be the subset
of $O'=\{o'_{1},...,o'_{n}\}$ containing all those $o'_{i}$ such
that $L_{j}$ intersects with $D(o'_{i},\frac{\delta'}{6})$, or equivalently,
with $B(o'_{i},\frac{\delta'}{6})$. That is,
\begin{equation}
\mathcal{C}_{j}:=\{o'_{i}\in O'\big|L_{j}\cap B(o'_{i},\frac{\delta'}{6})\neq\emptyset\}.\label{eq:def =00005Cmathcal=00007BC=00007D_j}
\end{equation}
Re-label the points in $\mathcal{C}_{j}$ as $o'_{(j,1)}$,..., $o'_{(j,n_{j})}$
in such a way that the corresponding longitude $w_{(j,1)}$, ... ,
$w_{(j,n_{j})}$ are counter clockwise on $S^{1}$. By the way the
circles $L_{j}$ are chosen, their diameters are strictly larger than
$\frac{\delta'}{3}$. Therefore, the intersection of $L_{j}$ with
each neighbourhood $D(o'_{(j,l)},\frac{\delta'}{6})$ (or equivalently,
$B(o'_{(j,l)},\frac{\delta'}{6})$) is an arc/interval $I_{(j,l)}$
centred at a point $q_{j,l}$. Here, $q_{j,l}$ is the intersection
of the longitude $R_{w_{(j,l)}}$ with the latitude. 
\begin{conclusion}
The distance between $o'_{(j,l)}$ and $q_{j,l}$ is less than $\frac{\delta'}{6}$.
Since $\{I_{(j,l)}\}$ covers $L_{j}$, the distance between each
consecutive points $q_{j,l}$ and $q_{j,l+1}$ is less than $\frac{\delta'}{3}$.
Here, $l\in\{1,...,n_{j}\}$ and $q_{j,l+1}$ refer to $q_{j,1}$
when $l=n_{j}$.
\end{conclusion}

For each $j=-k,...,k$, we set
\begin{equation}
\mathcal{V}_{j}:=\{q_{j,l}\big|,l=1,...,n_{j}\}.\label{eq:def=00005Cmathcal=00007BV=00007D_j}
\end{equation}

For $j<k-1$, we will triangulate the annulus $\mathcal{A}_{j,j+1}$
between $L_{j}$ and $L_{j+1}$ with the points in $\mathcal{V}_{j}$
and $\mathcal{V}_{j+1}$ as the vertices. T triangulate the discs
$D_{k}$ and $D_{-k}$ enclosed by respectively $L_{k}$ and $L_{-k}$,
we still need two more points lying in the discs. For $D_{k}$ we
pick $o'_{n}$. For $D_{-k}$, we choose a point $o'_{s}$ such that
$(0,0,-1)\in D(o'_{s},\frac{\delta'}{6})$. Since the distance from
$(0,0,-1)$ to the circle $L_{-k}$ is no less than $\frac{\delta'}{3}$,
we know that $o'_{s}$ lies in disc $D_{-k}$. Moreover, the neighbourhoods
$D(o'_{(-k,l)},\frac{\delta'}{6})$ do not cover the south pole $(0,0,-1)$
and then $o'_{s}$ is not in $\mathcal{C}_{-k}$ (or any other $\mathcal{C}_{j}$).
We take the set of vertices to be
\[
\mathcal{V}:=\big(\bigcup_{j}\mathcal{V}_{j}\big)\bigcup\{o'_{n},o'_{s}\}.
\]

\subsubsection{Connecting the vertices and constructing the triangulation}

Now we decribe the construction of the triangulation. The latitude
lines decompose the sphere into areas $\mathcal{A}_{j,j+1}$ between
$L_{j}$ and $L_{j+1}$ for $j\leq k-1$ and discs $D_{k}$ and $D_{-k}$
on $S^{2}$ respectively enclosed by $L_{k}$ and $L_{-k}$. The triangulation
of $S^{2}$ to be constructed will also triangulate each of these
areas. 

\textcolor{blue}{First of all, we connect $q_{j,l}$ with $q_{j,l+1}$
for each $j$ and $l$}. When $l=n_{j}$, $q_{j,l+1}$ refers to $q_{j,1}$.
We can simply take the closed arcs $[q_{j,l},q_{j,l+1}]$ to be the
edges. Note that the distance between $q_{j,l}$ and $q_{j,l+1}$
is less than $\frac{\delta'}{3}$, while the diameter of $L_{j}$
is larger than $\sqrt{2}\cdot\frac{\delta'}{3}$. Therefore, $[q_{j,l},q_{j,l+1}]$
is a minor arc on $L_{j}$.

Second, we connect\textcolor{blue}{{} $o_{n}$ to each $q_{k,l}$ in
$\mathcal{V}_{k}$ and $o'_{s}$ to each $q_{-k,l}$ in $\mathcal{V}_{-k}$.
}It indeed gives triangulations of $D_{k}$ and $D_{-k}$ since each
of the sets $\mathcal{V}_{k}$ and $\mathcal{V}_{-k}$ contains at
least $three$ points. We show this for $\mathcal{V}_{k}$ and it
works for $\mathcal{V}_{-k}$ in the same way. The radius of the circle
$L_{k}$ (in $\mathbb{R}^{3}$) is larger than $\frac{\sqrt{2}}{2}\cdot\frac{\delta'}{3}>\frac{\delta'}{6}$.
Meanwhile, the intersection of each $B(o'_{(k,l)},\frac{\delta'}{6})$
with the plane $\{z=h_{k}\}$ is a disc $D_{(k,l)}$ with radius no
more than $\frac{\delta'}{6}$. If $\mathcal{V}_{k}$ has only two
points, it means that $D_{(k,1)}$ and $D_{(k,2)}$ should cover $L_{k}$.
Then by symmetry, at least one of the discs should cover both the
ends of a diameter of $L_{k}$, which is impossible due to their sizes.
Moreover, the distance from $o_{n}$ to $L_{k}$ is less than $\frac{2\delta'}{3}$,
and by the triangle inequality $o'_{s}$ is no farther than $\frac{2\delta'}{3}+\frac{\delta'}{6}$
from each $q_{-k,l}$ since $(0,0,-1)\in D(o'_{s},\frac{\delta'}{6})$. 

Last, we complete the triangulation of the area $\mathcal{A}_{j,j+1}$
between $L_{j}$ and $L_{j+1}$ by connecting the points from $\mathcal{V}_{j}$
with those $\mathcal{V}_{j+1}$ in a proper way. For $q\in R_{w}$,
$q'\in R_{w'}$ and $q''\in R_{w''}$, we say that $q'$ is $\emph{before}$
$q''$ from $\emph{the right-hand side}$ of $q$ if and only if the
relations $w'=e^{i\theta'}w$ and $w''=e^{i\theta''}w$ hold for $0<\theta'<\theta''<2\pi$.
\textcolor{blue}{For $j=-k,...,k-1$, we connect each point $q_{j,l}$
to the first point in $\mathcal{V}_{j+1}$ from the right-hand side
of $q_{j,l}$, and $q_{j+1,l'}$ to the first point in $\mathcal{V}_{j}$
from the right-hand side of $q_{j+1,l'}$}\textcolor{purple}{. }We
explain in detail for $j\geq0$ that this gives a triangulation of
the area $\mathcal{A}_{j,j+1}$ and the distance between each pair
of connected points is less than $\frac{2\delta'}{3}$. The explanation
also works for the other cases in the same way. 

For $j\geq0$, the radius of $L_{j}$ is larger then $L_{j+1}$. Let
$q'_{j,l}$ be the point on $L_{j+1}$ with the same longitude with
$q_{(j,l)}$ . Observe that the distance between $q'_{j,l}$ and $q'_{j,l+1}$
is smaller than that between $q_{j,l}$ and $q_{j,l+1}$, and hence
is smaller than $\frac{\delta'}{3}$. Now we have two families of
points on $L_{j+1}$: $\{q_{j+1,l'}\}$ and $\{q'_{j,l}\}$. For convenience,
we color $q_{j+1,l'}$ in grey and $q'_{j,l}$ in red. The circle
$L_{j+1}$ can be decomposed into intervals each of which contains
points in a single color, and each point from $\{q_{j+1,l'}\}$ and
$\{q'_{j,l}\}$ belongs to one of these intervals. Suppose that $I^{r}$
is an interval contains points in $red$. Then it holds $I^{r}\subset(q_{j+1,l'},q_{j+1,l'+1})$
for some $l'\in\{1,...,n_{j+1}\}$. Since $(q_{j+1,l'},q_{j+1,l'+1})$
is a minor arc on $L_{j+1}$, each points $q'_{j,l}$ in $I^{r}$
is at a distance less than $\frac{\delta'}{3}$ from $q_{j+1,l'+1}$,
which is exacly the first point in $\mathcal{V}_{j+1}$ from the right-hand
side of $q_{j,l}$. Similarly, an interval $I^{g}$ containing only
grey points is contained in $(q'_{j,l},q'_{j,l+1})$ for some $l\in\{1,...,n_{j}\}$,
and then any point $q_{j+1,l'}$ in $I^{g}$ is no farther than $\frac{\delta'}{3}$
from $q'_{j,l+1}$, and $q_{j,l+1}$ is the first point in $\mathcal{V}_{j}$
from the right-hand side of $q_{j+1,l'}$. Since the distance between
$L_{j}$ and $L_{j+1}$ is $\frac{\delta'}{3}$, by the triangle inequality
we know that the distances between the connected pairs are indeed
less than $\frac{2\delta'}{3}$. 

\subsection{Conclusion}

We describe the sub-network $[p'_{ij}]_{0}$ obtained from the triangulation
constructed above. Take the subset of $\{o'_{1},...,o'_{n}\}$
\[
\mathcal{C}:=\big(\bigcup_{j}\mathcal{C}_{j}\big)\bigcup\{o'_{n},o'_{s}\},
\]
and then the points in $\mathcal{C}$ are in one-one correspondence
with the points in $\mathcal{V}$. (Recall that each point in $\mathcal{C}_{j}$
is re-labelled as $o'_{(j,l)}$ and corresponds to $q_{j,l}$ in $\mathcal{V}_{j}$,
and the points $o'_{n},\,o'_{s}$ simply correspond to themselves).
For $i,i'\in\{1,...,n\}$, if $o'_{i}$ and $o'_{i'}$ are in $\mathcal{C}$
and the corresponding points in $\mathcal{V}$ are connected in the
triangulation, we set $p'_{ii'}=1$, otherwise, set $p'_{ii'}=0$. 

We still need to verify that when $p'_{ii'}=1$, it holds $|o'_{i}-o'_{i'}|<\delta'$.
Suppose that $o'_{i},o'_{i'}\in\big(\bigcup_{j}\mathcal{C}_{j}\big)$
and correspond to $q_{i,l}$ and $q_{i',l'}$ in $\big(\bigcup_{j}\mathcal{V}_{j}\big)$,
respectively. Then when they are connected, we have $|q_{i,l}-q_{i',l'}|<\frac{2\delta'}{3}$.
The distances between $q_{j,l}$ and $o'_{(j,l)}$ are less than $\frac{\delta'}{6}$
for all $(j,l)$, and hence it holds 
\[
|o'_{i}-o'_{i'}|<|o'_{(i,l)}-q_{i,l}|+|q_{i,l}-q_{i',l'}|+|q_{i',l'}-o'_{(i',l')}|<\delta'.
\]
If $o'_{i}=o'_{n}$ and $o'_{i'}=o'_{(k,l)}$, it holds $|o'_{i}-q_{k,l}|<\frac{2\delta'}{3}+\frac{\delta'}{6}$
and then again we get $|o'_{i}-o'_{i'}|<\delta'$ from the triangle
inequality. The same works for the case with $o'_{i}=o'_{s}$ and
$o'_{i'}=o'_{(-k,l)}$.

Combining Propositions \ref{prop:main result} and \ref{prop:natural condition and triangulation},
we conclude that
\begin{thm}
\label{thm:network-variation under natural-condition}Suppose that
the initial positions $\{o_{1},...,o_{n}\}$ of the robots constitute
a $\frac{\delta}{6}$-net on $S^{2}$, and any two robots will set
up a message channel if the distance between them is smaller than
$\delta$. If after a movement $\mathcal{R}$ the robots land on $\mathfrak{P}$
at the moment $\bar{t}>0$, then $[p_{ij}]_{t}\neq[p_{ij}]_{0}$ for some $t\in[0,\bar{t}]$.
\end{thm}

\bibliographystyle{plain}
\bibliography{ControlTheory.bib}
 
\end{document}